\documentclass[preprint,a4paper,12pt]{elsarticle}

\usepackage{lineno}

\usepackage{amsmath}
\usepackage{amsfonts}
\usepackage{amssymb}
\usepackage{amsthm}

\usepackage{graphicx}
\usepackage{subcaption}
\usepackage{multirow} 
\usepackage{float} 

\theoremstyle{definition}

\usepackage{xcolor}
\usepackage[colorlinks=true, allcolors=blue]{hyperref}

\newcommand{\abs}[1]{\left|#1\right|}
\newcommand{\norm}[1]{\left\lVert#1\right\rVert}

\journal{Physica D: Nonlinear Phenomena}

\begin{document}
	
\begin{frontmatter}
	
	
	
	\title{Physics-informed neural networks for high-dimensional solutions and snaking bifurcations in nonlinear lattices}
	
	
	
	\author[inst1]{Muhammad Luthfi Shahab}
	
	\author[inst2]{Fidya Almira Suheri}
	\author[inst2]{Rudy Kusdiantara}
	\author[inst1]{Hadi Susanto}
	
	\affiliation[inst1]{
		organization={Department of Mathematics, Khalifa University of Science \& Technology},
		city={Abu Dhabi},
		postcode={PO Box 127788}, 
		country={United Arab Emirates}
	}
	
	\affiliation[inst2]{
		organization={Industrial and Financial Mathematics Research Group, Institut Teknologi Bandung},
		city={Bandung},
		postcode={40132}, 
		country={Indonesia}
	}
	
	
\begin{abstract}
This paper introduces a framework based on physics-informed neural networks (PINNs) for addressing key challenges in nonlinear lattices, including solution approximation, bifurcation diagram construction, and linear stability analysis. We first employ PINNs to approximate solutions of nonlinear systems arising from lattice models, using the Levenberg–Marquardt algorithm to optimize network weights for greater accuracy. To enhance computational efficiency in high-dimensional settings, we integrate a stochastic sampling strategy. We then extend the method by coupling PINNs with a continuation approach to compute snaking bifurcation diagrams, incorporating an auxiliary equation to effectively track successive solution branches. For linear stability analysis, we adapt PINNs to compute eigenvectors, introducing output constraints to enforce positivity, in line with Sturm–Liouville theory. Numerical experiments are conducted on the discrete Allen–Cahn equation with cubic and quintic nonlinearities in one to five spatial dimensions. The results demonstrate that the proposed approach achieves accuracy comparable to, or better than, traditional numerical methods, especially in high-dimensional regimes where computational resources are a limiting factor. These findings highlight the potential of neural networks as scalable and efficient tools for the study of complex nonlinear lattice systems.
\end{abstract}

	
	
	
\end{frontmatter}


\section{Introduction}

Neural networks have emerged as powerful tools for solving partial differential equations (PDEs) \cite{cuomo2022scientific}, largely due to their capacity to approximate complex functions with high accuracy \cite{hornik1989multilayer}. Foundational studies by Dissanayake and Phan-Thien \cite{dissanayake1994neural}, and Lagaris et al.\ \cite{lagaris1998artificial}, demonstrated that PDEs can be reformulated as optimization problems and solved using neural network architectures, marking a significant shift in computational mathematics. More recently, physics-informed neural networks (PINNs) \cite{raissi2019physics} have advanced the field by incorporating governing physical laws directly into the training process, enabling the resolution of both forward and inverse problems across a range of domains, including fluid dynamics \cite{cai2021physics} and optics \cite{chen2020physics}. A key strength of neural networks lies in their ability to overcome the curse of dimensionality. Unlike traditional numerical methods, which become computationally infeasible in high-dimensional settings, neural networks offer scalable and efficient alternatives through optimization-based formulations \cite{han2018solving, putri2024deep}, thereby enabling the solution of previously intractable problems.

Extending these capabilities, recent studies have applied neural network methodologies to nonlinear lattices—discrete systems that arise in many physical contexts. Zhu et al.\ \cite{zhu2022neural} employed symmetry-preserving neural architectures to model nonlinear dynamical lattices, focusing on the Ablowitz–Ladik system. Lin and Chen \cite{lin2024pseudo} introduced a pseudo-grid-based, physics-informed convolutional-recurrent network for integrable lattice equations, while Geng et al.\ \cite{geng2024separable} developed a separable graph Hamiltonian network to capture lattice dynamics via graph-based deep learning. Zhou et al.\ \cite{zhou2024symmetric} integrated symmetric difference data enhancement into a PINN framework to improve accuracy for discrete nonlinear equations.

In parallel, alternative paradigms have broadened the scope of neural approaches. Saqlain et al.\ \cite{saqlain2023discovering} proposed PINNs for discovering governing equations in discrete systems. Opala et al.\ \cite{opala2019neuromorphic} applied reservoir computing to the complex Ginzburg–Landau lattice model, successfully capturing a wide range of nonlinear dynamics. Stokes et al.\ \cite{stokes2020phases} introduced first-quantized deep neural networks for analyzing strongly coupled fermionic lattice systems, further expanding the application of neural networks in quantum many-body physics. Together, these efforts underscore the growing versatility of deep learning techniques in tackling the complexities of nonlinear lattice models.

Nonlinear lattices play a central role in modeling phenomena in optics, condensed matter, and biological systems \cite{kartashov2011solitons}. However, they pose substantial analytical and computational challenges, particularly in bifurcation analysis and stability characterization—both essential for understanding qualitative system behavior. These difficulties are compounded by intricate bifurcation structures and bistability, common features in discrete settings \cite{taylor2010snaking, kusdiantara2017homoclinic}, and demand precise, scalable computational strategies.

This study addresses a key gap in the application of PINNs to nonlinear lattice systems, with emphasis on high-dimensional cases, bifurcation structures, and linear stability analysis. While PINNs have been successfully applied to related PDE-based tasks \cite{shahab2024neural, shahab2025corrigendum, fabiani2021numerical, galaris2022numerical}, their use in discrete lattice contexts remains limited. Here, we extend PINN methodologies to these systems, particularly in the context of bifurcation tracking and stability diagnostics.

We also consider high-dimensional problems, where conventional solvers often fail due to the curse of dimensionality. PINNs offer an efficient alternative capable of capturing complex dynamics in such settings. Our analysis centers on the discrete Allen–Cahn equation with cubic and quintic nonlinearities \cite{taylor2010snaking}, a benchmark model known for its rich bifurcation landscape, including snaking patterns with multiple turning points. This system provides a rigorous testbed for evaluating the effectiveness of neural-network-based solvers in nonlinear lattice dynamics.

The remainder of this paper is organized as follows. Section~\ref{sec2} introduces the discrete Allen–Cahn equation. Section~\ref{sec3} discusses the pseudo-arclength continuation method. Section~\ref{sec4} presents the proposed PINN framework, including its integration with pseudo-arclength continuation, a stochastic approach for high-dimensional problems, and a linear stability analysis module. Section~\ref{sec5} provides numerical results from one-dimensional to five-dimensional cases, demonstrating the versatility and performance of the proposed approach. Finally, Section~\ref{sec6} summarizes the key findings and outlines directions for future research.

\section{Discrete Allen--Cahn Equation}
\label{sec2}

We consider the discrete Allen--Cahn equation with cubic and quintic nonlinearities and a parameter $\mu$. This model arises from the continuous PDE
\begin{equation}
	\frac{\partial u}{\partial t} = \mu u + \Delta u + 2u^3 - u^5,
\end{equation}
by discretizing the spatial domain. Specifically, the continuous variable $u$ is replaced by its discrete counterpart $u_i$ defined on a one-dimensional lattice, $i = 1, \dots, n$, and the Laplacian $\Delta u$ is approximated using the finite difference operator $c \Delta_1 u_i = c (u_{i+1} - 2u_i + u_{i-1})$, where $c > 0$ is a fixed coupling parameter. This yields the one-dimensional discrete Allen--Cahn equation:
\begin{equation} \label{allencahn1}
	\frac{du_i}{dt} = \mu u_i + c\Delta_1 u_i + 2u_i^3 - u_i^5, \quad i \in \{2, \dots, n-1\},
\end{equation}
with boundary conditions $u_1 = u_n = 0$.

The two-dimensional version is defined on a lattice $(i,j)$, where $i,j \in \{1, \dots, n\}$. The equation becomes
\begin{equation} \label{allencahn2}
	\frac{du_{i,j}}{dt} = \mu u_{i,j} + c\Delta_2 u_{i,j} + 2u_{i,j}^3 - u_{i,j}^5, \quad i,j \in \{2, \dots, n-1\},
\end{equation}
where
\begin{equation}
	c \Delta_2 u_{i,j} = c(u_{i+1,j} + u_{i-1,j} + u_{i,j+1} + u_{i,j-1} - 4u_{i,j}),
\end{equation}
and boundary conditions are $u_{i,j} = 0$ for $i \in \{1,n\}$ or $j \in \{1,n\}$. Higher-dimensional versions (three to five dimensions) can be defined analogously by extending the discrete Laplacian accordingly.

In all cases, we are interested in the steady-state solutions obtained by setting the time derivatives to zero. This reduces the system to a set of nonlinear algebraic equations. For one- and two-dimensional cases, we additionally analyze the bifurcation structure and linear stability of the steady states. The discrete Allen--Cahn model is particularly well-suited for this purpose due to its characteristic snaking bifurcations with multiple turning points \cite{taylor2010snaking}, making it an effective benchmark for continuation methods.

The steady states of interest are symmetric about the midpoint of the lattice \cite{taylor2010snaking}. To capture this, the lattice size is commonly chosen as $n = 2m-1$ for site-centered solutions, or $n = 2m$ for bond-centered solutions. These configurations allow for a detailed comparison of solution types and their respective bifurcation characteristics.

\subsection{Eigenvalue Problem for Linear Stability}
\label{sec2a}

To determine the linear stability of a steady-state solution $u$ to Eq.~\eqref{allencahn1}, we consider a small perturbation of the form
\begin{equation}
	\Tilde{u}_i(t) = u_i + \varepsilon w_i(t), \qquad \abs{\varepsilon} \ll 1.
\end{equation}
Substituting this into Eq.~\eqref{allencahn1} and linearizing in $\varepsilon$ yields the equation for the perturbation:
\begin{equation}
	\frac{dw_i}{dt} = \mu w_i + c\Delta_1 w_i + 6 u_i^2 w_i - 5 u_i^4 w_i, \quad i \in \{2, \dots, n-1\}.
\end{equation}
Assuming a separable solution $w_i(t) = e^{\lambda t} v_i$, we arrive at the eigenvalue problem:
\begin{equation} \label{allencahn_eigenvalue}
	\lambda v_i = \mu v_i + c\Delta_1 v_i + 6 u_i^2 v_i - 5 u_i^4 v_i, \quad i \in \{2, \dots, n-1\},
\end{equation}
subject to boundary conditions $v_1 = v_n = 0$, where $\lambda$ is the eigenvalue and $v = (v_i)$ is the corresponding eigenvector.

A steady-state solution is linearly stable if the largest eigenvalue satisfies $\max(\lambda) < 0$; otherwise, it is unstable. The formulation for higher dimensions follows analogously. Notably, Eq.~\eqref{allencahn_eigenvalue} can be viewed as a discrete Sturm--Liouville problem. 
To determine $\max(\lambda)$, we identify the critical eigenvalue associated with an eigenvector that has no internal sign changes.

\section{Pseudo-arclength Continuation}
\label{sec3}

In parameter-dependent nonlinear systems, direct continuation by varying the parameter may fail to trace the full bifurcation structure, especially in the presence of turning points. At such points, the solution curve folds back on itself, and the parameter is no longer a monotonic or uniquely defining quantity for the solution. Standard continuation methods often fail here due to the breakdown of the one-to-one correspondence between the parameter and the solution.

Pseudo-arclength continuation overcomes this limitation by reparameterizing the solution curve using a pseudo-arclength measure, allowing the method to follow the solution branch through turning points. Instead of advancing the parameter directly, the method controls the step direction along the curve itself via an auxiliary constraint.

Consider a system of nonlinear equations $F(u,\mu)=0$, where $u$ is the state variable and $\mu$ is a bifurcation parameter. For simplicity, denote a solution as $(u^{[k]}, \mu^{[k]})$, where $u^{[k]}$ satisfies $F(u, \mu^{[k]}) = 0$. To construct the bifurcation diagram, we begin with two known solutions, $(u^{[1]}, \mu^{[1]})$ and $(u^{[2]}, \mu^{[2]})$, and generate subsequent solutions using a predictor-corrector scheme.

The predictor step uses the extrapolation:
\begin{equation} \label{predictor}
	(u^{[k]}, \mu^{[k]}) = 2(u^{[k-1]}, \mu^{[k-1]}) - (u^{[k-2]}, \mu^{[k-2]}),
\end{equation}
which provides an initial guess for the $k$-th solution.
To enable correction and enforce continuation along the solution branch, we introduce a constraint equation of the form:
\begin{equation} \label{arclength1}
	\alpha \left( \sqrt{ \left[\beta_1 \left(\norm{u^{[k]}} - (\norm{u^{[k-1]}} + \gamma) \right) \right]^2 + \left[\beta_2 (\mu^{[k]} - \mu^{[k-1]}) \right]^2 } - \delta \right) = 0.
\end{equation}
This auxiliary condition serves two main purposes, depending on the choice of parameters:
\begin{itemize}
	\item \textbf{Norm-based continuation}: When $\beta_2 = \delta = 0$, $\beta_1 = 1$, and $\gamma \ne 0$, the continuation enforces $\norm{u^{[k]}} = \norm{u^{[k-1]}} + \gamma$. This approach is effective when each solution norm corresponds to a unique solution, as seen in problems such as the Bratu equation \cite{joseph1973quasilinear}, the Burgers equation \cite{allen2013numerical, fabiani2021numerical}, and nonlinear boundary value problems \cite{graef2003three, liao2012homotopy}. We adopt this formulation for the one-dimensional discrete Allen--Cahn equation.

	\item \textbf{Arclength-based continuation}: When $\beta_1, \beta_2 \ne 0$ and $\gamma = 0$, the constraint enforces a fixed pseudo-arclength $\delta$ between successive steps in the $(\norm{u}, \mu)$ plane. This formulation is better suited for general bifurcation tracking when the norm and parameter operate on different scales. We apply this approach for the two-dimensional discrete Allen--Cahn equation.
\end{itemize}
Here, $\alpha$ is a weighting factor balancing the continuation constraint with the original system $F(u, \mu) = 0$. When $\delta = 1$, the square root in Eq.~\eqref{arclength1} may be omitted for computational simplicity. The effects of tuning parameters $\alpha$, $\gamma$, and $(\beta_1, \beta_2)$ are discussed in Appendices~\ref{appendix_alpha_gamma} and~\ref{appendix_beta}.

The continuation proceeds iteratively until a stopping condition is met (e.g., reaching $k = k_{\max}$ or covering the desired range of solutions). For improved convergence, the initial solutions $(u^{[1]}, \mu^{[1]})$ and $(u^{[2]}, \mu^{[2]})$ should approximately satisfy Eq.~\eqref{arclength1}.

Each correction step requires solving the augmented system consisting of $F(u, \mu) = 0$ and Eq.~\eqref{arclength1}. Newton's method can be applied, with iterations given by \cite{chong2013introduction}:
\begin{equation}
	u^{(k+1)} = u^{(k)} - J(u^{(k)})^{-1} F(u^{(k)}),
\end{equation}
where $J(u^{(k)})$ is the Jacobian of $F$ evaluated at $u^{(k)}$. The dominant computational cost lies in solving the linear system involving $J(u^{(k)})$.

In this work, we employ the Levenberg--Marquardt algorithm \cite{levenberg1944method, marquardt1963algorithm}, a regularized variant of Newton's method, given by:
\begin{equation}
	u^{(k+1)} = u^{(k)} - (J^T J + \lambda^{(k)} I)^{-1} J^T F(u^{(k)}),
\end{equation}
where $\lambda^{(k)}$ is a damping parameter ensuring that each iteration proceeds in a descent direction. When solving the combined system $F(u, \mu) = 0$ and Eq.~\eqref{arclength1}, the parameter $\mu$ is treated as an additional variable and updated alongside $u$.

We use $u^{(k)}$ to denote the $k$-th iterate produced by the solver, and $u^{[k]}$ to represent the $k$-th converged solution that satisfies $F(u, \mu) = 0$. This notation distinguishes from $u^k$, which denotes exponentiation.

\section{Proposed Methods}
\label{sec4}

\subsection{PINNs for Nonlinear Lattices}
\label{sec_4_1}

Instead of solving the nonlinear system directly using traditional numerical solvers, this study proposes the use of physics-informed neural networks (PINNs) to approximate solutions. While PINNs have been extensively applied to continuous partial differential equations (PDEs), their adaptation to discrete nonlinear lattice problems remains comparatively underexplored.

We begin with a general nonlinear system of equations of the form
\begin{equation} \label{system1}
	f_{\mathbf{i}}(u, \mu) = 0, \qquad \mathbf{i} \in A,
\end{equation}
where $u = \{ u_{\mathbf{i}} \mid \mathbf{i} \in A \}$ denotes the set of unknowns, $\mu$ is a fixed parameter, and $A$ is the index set corresponding to the discrete lattice points. The multi-index $\mathbf{i}$ may represent one or more spatial indices, such as $(i_1, i_2)$ or $(i, j)$, depending on the system dimension.

Unlike continuous PDE problems where the network input corresponds to spatial or temporal coordinates, nonlinear lattice systems lack explicit spatial variables. To overcome this, we reinterpret the lattice index $\mathbf{i}$ as an input to the network and denote the network output as $u(\mathbf{i})$.

Given this formulation, a neural network is constructed to approximate the solution as $u(\mathbf{i}, W)$, where $W$ represents the weights and biases of the network. The goal is to find $W$ such that
\begin{equation} \label{nonlinear_system_NN}
	f_{\mathbf{i}}(u(\mathbf{i}, W), \mu) = 0, \qquad \mathbf{i} \in A.
\end{equation}
This system can be interpreted as a root-finding problem in the space of network parameters. To measure the discrepancy, we define the mean squared error (MSE) as:
\begin{equation} \label{MSE_NN}
	\text{MSE} = \frac{1}{|A|} \sum_{\mathbf{i} \in A} \left[ f_{\mathbf{i}}(u(\mathbf{i}, W), \mu) \right]^2,
\end{equation}
where $|A|$ denotes the number of lattice points in the system.

To determine the optimal weights $W$, two primary strategies can be employed: (1) minimize the loss function in Eq.~\eqref{MSE_NN} using optimization algorithms such as L-BFGS \cite{liu1989limited} or Adam \cite{kingma2014adam}, or (2) solve the system in Eq.~\eqref{nonlinear_system_NN} directly using a nonlinear solver. In this study, we adopt the second approach and apply the Levenberg--Marquardt algorithm \cite{levenberg1944method, marquardt1963algorithm}, which offers improved accuracy and convergence for small- to medium-sized networks.

The neural networks used in this work are shallow, comprising two hidden layers. Let $W_1$, $B_1$, $W_2$, $B_2$, $W_3$, and $B_3$ represent the weights and biases of the network. Then the output can be expressed as:
\begin{equation} \label{NN_function}
	u(\mathbf{i}, W) = \sigma\left( \sigma\left( \mathbf{i} W_1 + B_1 \right) W_2 + B_2 \right) W_3 + B_3,
\end{equation}
where $\sigma$ is the activation function. Following the success of previous work in nonlinear bifurcation problems \cite{shahab2024neural}, we adopt the Gaussian activation function $\sigma(x) = e^{-x^2/2}$ throughout this study.

\subsection{Input Manipulation}

In our implementation, we found that transforming the input index $\mathbf{i}$ before feeding it into the neural network significantly improves convergence. Experimental results indicate that mapping discrete indices to a smaller, normalized domain enhances the training efficiency and stability of PINNs. This improvement is motivated by two key observations:

\begin{itemize}
	\item Most activation functions perform optimally when their inputs lie within a bounded range. Large input indices can lead to extreme values when propagated through network layers, potentially causing vanishing gradients. Normalizing the input ensures that the network operates within the sensitive regions of its activation functions, which accelerates training and improves accuracy.

	\item PINNs are known to struggle with problems defined over large physical or index domains \cite{moseley2023finite, jagtap2020extended, heinlein2021combining}. These difficulties stem from highly non-convex loss landscapes and gradient pathologies. Scaling the input domain helps mitigate these issues and supports more robust optimization.
\end{itemize}

To formalize this procedure, consider a two-dimensional nonlinear lattice where $\mathbf{i} = (i_1, i_2)$ or $(i, j)$ and $i_1, i_2 \in \{1, \dots, n\}$. We apply the transformation:
\begin{equation} \label{eq_first_manipulation}
	\bar{i}_1 = (i_1 - 1)h, \qquad \bar{i}_2 = (i_2 - 1)h,
\end{equation}
where $h = 1 / (n - 1)$, so that the original index range $[1, n]$ is mapped to the normalized interval $[0, 1]$. The PINN then approximates the solution as $u(\mathbf{i}, W)$ using the transformed inputs.

Further refinements are possible when specific structural properties of the solution are known. For example, solutions of the discrete Allen--Cahn equation exhibit reflection symmetry about the midpoint of the lattice \cite{taylor2010snaking}. To exploit this, we propose an additional symmetry-aware transformation:
\begin{equation} \label{eq_second_manipulation}
	\begin{split}
		\hat{i}_1 = 0.5 - \abs{\bar{i}_1 - 0.5}, & \qquad \hat{i}_2 = 0.5 - \abs{\bar{i}_2 - 0.5}, \\
		(\tilde{i}_1, \tilde{i}_2) & = \text{sort}(\hat{i}_1, \hat{i}_2).
	\end{split}
\end{equation}
Here, $0.5$ represents the midpoint of the normalized domain. This transformation maps symmetrically equivalent inputs—those corresponding to equal solution values—to a canonical form within the region $0 \le \tilde{i}_1 \le \tilde{i}_2 \le 0.5$. As a result, it reduces the complexity of the optimization process, leading to faster convergence during training.

This preprocessing step incurs no additional cost during inference, as the transformations can be computed and cached before training begins. Further details and illustrations of this approach are provided in \cite{shahab2025finite}.

\subsection{Pseudo-arclength Continuation}
\label{sec_4_2}

When constructing bifurcation diagrams using pseudo-arclength continuation, the objective is to compute a sequence of solution-parameter pairs $(u^{[k]}, \mu^{[k]})$. However, since the solution $u$ is represented implicitly by a PINN through its weights $W$, we instead compute a sequence $(W^{[k]}, \mu^{[k]})$, where each $W^{[k]}$ defines an approximate solution $u(\mathbf{i}, W^{[k]})$.

To integrate pseudo-arclength continuation within the PINN framework, we first determine two initial solutions, $(W^{[1]}, \mu^{[1]})$ and $(W^{[2]}, \mu^{[2]})$, that satisfy the nonlinear system in Eq.~\eqref{nonlinear_system_NN}. For each subsequent step $k \geq 3$, the parameter $\mu^{[k]}$ is treated as an additional trainable variable alongside the neural network weights. The predictor step is initialized by extrapolation:
\begin{equation} \label{predictor_NN}
	(W^{[k]}, \mu^{[k]}) = 2(W^{[k-1]}, \mu^{[k-1]}) - (W^{[k-2]}, \mu^{[k-2]}).
\end{equation}
Using the predicted pair $(W^{[k]}, \mu^{[k]})$ as an initial guess, we retrain the network to satisfy both the nonlinear system and the continuation constraint:
\begin{equation} \label{arclength_NN}
	\begin{split}
		f_{\mathbf{i}}(u(\mathbf{i}, W^{[k]}), \mu^{[k]}) = 0, \qquad \mathbf{i} \in A, \\
		\alpha \left( \sqrt{ \left[ \beta_1 \left( \|u(\mathbf{i}, W^{[k]})\| - (\|u(\mathbf{i}, W^{[k-1]})\| + \gamma) \right) \right]^2 + \left[ \beta_2 (\mu^{[k]} - \mu^{[k-1]}) \right]^2 } - \delta \right) = 0.
	\end{split}
\end{equation}

Once the updated pair $(W^{[k]}, \mu^{[k]})$ is obtained, we evaluate $u(\mathbf{i}, W^{[k]})$ to plot the corresponding point on the bifurcation diagram. This process is iterated to trace the solution branch across turning points and parameter regimes of interest.

While the overall framework is conceptually straightforward, its implementation poses several technical challenges. These include:
\begin{itemize}
	\item Treating $\mu$ as a trainable parameter integrated into the network architecture.
	\item Initializing network weights via Eq.~\eqref{predictor_NN}, which is not standard in conventional machine learning workflows.
	\item Solving the nonlinear system using the Levenberg--Marquardt algorithm, as opposed to minimizing a typical loss function.
\end{itemize}

These requirements render standard machine learning libraries such as Scikit-learn \cite{pedregosa2011scikit} or high-level TensorFlow APIs \cite{abadi2016tensorflow} unsuitable. While functional (low-level) TensorFlow could support this approach, doing so would involve complex, non-standard manipulation of weights and control logic to implement Eqs.~\eqref{predictor_NN} and \eqref{arclength_NN} directly. For these reasons, we developed a custom implementation in MATLAB, which provides the flexibility necessary to accommodate the algorithmic structure and optimization routine specific to our method.

\subsection{Stochastic Newton's Method}
\label{stochastic_newton}

As discussed earlier, the PINN is trained to satisfy a nonlinear system using the Levenberg--Marquardt algorithm \cite{levenberg1944method, marquardt1963algorithm}, which is based on Newton’s method. In this subsection, we introduce a simple yet effective modification inspired by the mini-batching technique in stochastic gradient descent (SGD). Rather than using the full system of equations at every iteration, we operate on randomly selected subsets of equations, significantly reducing computational cost. This stochastic strategy is particularly applied to the five-dimensional discrete Allen--Cahn equation.

In the standard Newton's method, the network weights $W$ are iteratively updated as:
\begin{equation}
	W^{(k+1)} = W^{(k)} - J(W^{(k)})^{-1} F(W^{(k)}),
\end{equation}
where $J(W^{(k)})$ is the Jacobian of the nonlinear system $F$ evaluated at $W^{(k)}$. However, in high-dimensional settings, the computation and inversion of the full Jacobian become prohibitively expensive.

To address this, we adopt a stochastic approach. At each iteration, we randomly select a subset $S_k$ of the system equations and compute the corresponding Jacobian and residuals. This leads to the update rule:
\begin{equation}
	W^{(k+1)} = W^{(k)} - J_{S_k}(W^{(k)})^{-1} F_{S_k}(W^{(k)}),
\end{equation}
where $F_{S_k}$ and $J_{S_k}$ denote the restriction of $F$ and $J$ to the subset $S_k$. Although the approximation to the full Jacobian is less accurate, the reduced computational load per iteration enables faster convergence in practice.

This modification integrates naturally with neural networks since all weights contribute to computing $F_{S_k}$ and $J_{S_k}$. Similar ideas have been studied in the context of optimization \cite{byrd2016stochastic, kovalev2019stochastic}, though most previous work focuses on loss minimization rather than solving nonlinear systems directly.

When using the Levenberg--Marquardt algorithm, the standard update becomes:
\begin{equation}
	W^{(k+1)} = W^{(k)} - (J^T J + \lambda^{(k)} I)^{-1} J^T F(W^{(k)}),
\end{equation}
which we adapt to the stochastic setting as:
\begin{equation}
	W^{(k+1)} = W^{(k)} - \left( J_{S_k}^T J_{S_k} + \lambda^{(k)} I \right)^{-1} J_{S_k}^T F_{S_k}(W^{(k)}).
\end{equation}

It is important to note that this stochastic method is particularly suited to neural-network-based formulations. In traditional finite-difference or finite-element contexts, functions in $F$ are typically local, depending only on adjacent variables. Therefore, randomly sampling equations may fail to update all variables, limiting convergence. Neural networks, by contrast, maintain global parameter coupling, making the stochastic method viable for high-dimensional problems.

\subsection{Computing the Nonlinear System and Treating Boundary Points}

\begin{figure}[t]
	\centering
	\includegraphics[width=\textwidth]{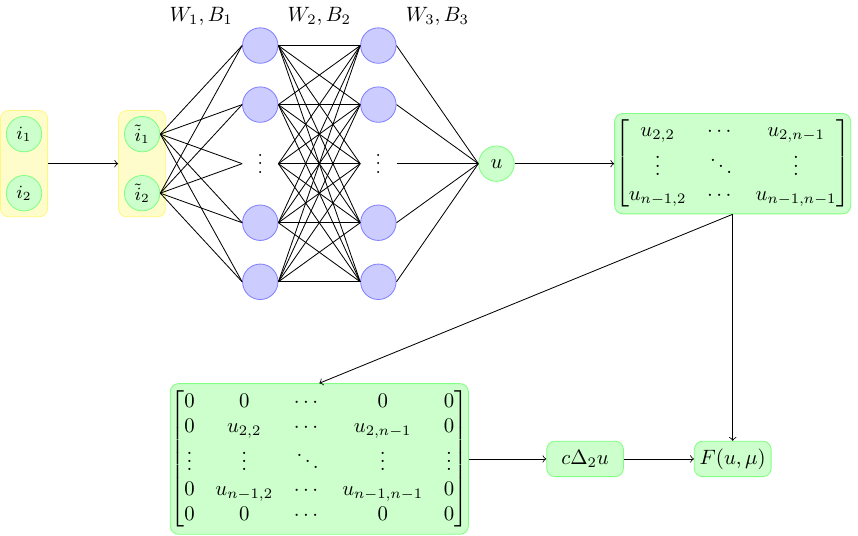}
	\caption{Schematic of the PINN architecture for solving nonlinear lattice systems. Although the diagram illustrates a two-dimensional lattice, the method generalizes to higher dimensions by extending the input space accordingly.}
	\label{fig_PINNs}
\end{figure}

For one- to four-dimensional problems, we evaluate the full set of nonlinear equations in each iteration of the Levenberg--Marquardt algorithm. A $d$-dimensional matrix $u$ is first constructed to represent the solution domain, where $d$ is the number of spatial dimensions. Boundary conditions are enforced by explicitly setting $u_{\mathbf{i}} = 0$ at all boundary indices, eliminating the need for neural network evaluation at those locations. For interior points $\mathbf{i} \in A$, the PINN is used to compute $u_{\mathbf{i}} = u(\mathbf{i}, W)$. This approach allows the nonlinear system to be efficiently evaluated using matrix operations.

In the five-dimensional case, due to memory and computational constraints, we employ the stochastic Newton's method described above. At each iteration, only a subset $S_k$ of the nonlinear equations is evaluated. Consequently, we avoid forming a full five-dimensional matrix. Instead, we construct an input vector corresponding to the points in $S_k$, along with 10 additional vectors representing the neighboring lattice sites. The PINN then computes the solution values $u(\mathbf{i}, W)$ for all sampled and adjacent points.

Some of these indices may lie on the boundary. To enforce zero boundary conditions in this case, we apply a multiplicative mask to the neural network output using:
\begin{equation}
	u(\mathbf{i}, W) \gets u(\mathbf{i}, W) \cdot \prod_{k=1}^5 \sin(\pi i_k),
\end{equation}
for $\mathbf{i} = (i_1, i_2, i_3, i_4, i_5)$. This transformation ensures that the output vanishes on the boundary while preserving flexibility in the interior. The nonlinear system is evaluated using vectorized operations for scalability and computational efficiency.

Figure~\ref{fig_PINNs} illustrates the general structure of the PINN used in our framework. While the diagram focuses on a two-dimensional lattice for clarity, the extension to higher dimensions is straightforward: the number of input coordinates increases, and the internal matrix representation is adapted accordingly.

\subsection{Largest Eigenvalue}
\label{sec_4_3}

Consider the general eigenvalue problem
\begin{equation}
	H(v, u, \mu) = \lambda v,
\end{equation}
where $u$ and $\mu$ are fixed, $\lambda$ is an eigenvalue, and $v$ is the corresponding eigenvector. This equation can be reformulated as a nonlinear system:
\begin{equation}
	G(v, \lambda, u, \mu) = H(v, u, \mu) - \lambda v = 0.
\end{equation}

To solve this system, we adopt a PINN-based approach similar to that described in Section~\ref{sec_4_1}. Our objective is to compute the \emph{largest} eigenvalue, which is critical for linear stability analysis as noted in Section~\ref{sec2a}. According to Sturm--Liouville theory, the dominant eigenvalue is associated with an eigenvector that does not change sign. To enforce this non-negativity, we apply an absolute value transformation to the neural network output. Using the same notation for weights and biases, the neural network is defined as:
\begin{equation} \label{NN_eigenproblem}
	v(\mathbf{i}, W) = \abs{\sigma\left(\sigma\left(\mathbf{i} W_1 + B_1\right) W_2 + B_2\right) W_3 + B_3},
\end{equation}
where $\sigma$ denotes the activation function.

If the system $G$ consists of multiple component equations $g_{\mathbf{i}}$, indexed over a set $\mathbf{i} \in A$, the PINN is trained to satisfy:
\begin{equation} \label{NN_eigenvalue}
	\begin{split}
		g_{\mathbf{i}}(v(\mathbf{i}, W), \lambda, u, \mu) &= 0, \qquad \mathbf{i} \in A, \\
		\|v\|_2 &= 1,
	\end{split}
\end{equation}
where the normalization constraint $\|v\|_2 = 1$ ensures uniqueness of the eigenvector. In this setup, $\lambda$ is treated as an additional trainable parameter, jointly optimized with the network weights $W$.

Upon convergence, the solution $\lambda$ obtained from Eq.~\eqref{NN_eigenvalue} corresponds to the largest eigenvalue, consistent with the Sturm--Liouville theory. This PINN-based formulation offers a general and flexible framework for computing dominant eigenvalues in high-dimensional, structured systems.

\section{Experimental Results and Discussion}
\label{sec5}

We use the Levenberg--Marquardt algorithm \cite{levenberg1944method, marquardt1963algorithm}, both for direct solution of the nonlinear system and via the PINN-based formulations described in Section~\ref{sec4}. Solutions obtained through direct numerical methods serve as ground truth for benchmarking the accuracy of the PINN approximations. All simulations were performed using MATLAB R2024a on a laptop equipped with an Intel Core i7-1185G7 processor.

\subsection{One-dimensional Discrete Allen--Cahn Equation}

We begin with the one-dimensional discrete Allen--Cahn equation, defined in Eq.~\eqref{allencahn1}. It is well established that this equation exhibits characteristic snaking bifurcation structures. For our experiments, we set $m=10$, resulting in lattice sizes of $n=19$ and $n=20$ for the site- and bond-centered configurations, respectively. Excluding boundary points, the Jacobian matrices for direct solution have dimensions $17 \times 17$ and $18 \times 18$.

For the PINN approach, we use a network architecture with one input, two hidden layers of four neurons each, and one output. This yields a total of 33 trainable parameters, including biases. The Jacobian matrices for this configuration expand to $17 \times 33$ and $18 \times 33$ for the site- and bond-centered cases, respectively. Reducing the network size adversely affects accuracy. Notably, the relative increase in Jacobian size occurs only in low-dimensional settings with small lattice sizes; in higher-dimensional problems, the PINN Jacobians are often more compact relative to the original system size.

Bifurcation diagrams are constructed using the norm-based continuation with the following parameters: $\alpha = 10$, $\beta_1 = 1$, $\beta_2 = 0$, $\gamma = 10/1000$, and $\delta = 0$. The solution norm is defined as
\begin{equation}
	\norm{u} = \frac{1}{1+\sqrt{1+\mu}} \sum_{i=1}^{n} u_i^2,
\end{equation}
and diagrams are computed up to $\norm{u} = 10$ for both solution types. Sensitivity to $\alpha$ and $\gamma$ is explored in Appendix~\ref{appendix_alpha_gamma}.

\begin{figure}
	\centering
	\begin{subfigure}[b]{0.49\textwidth}
		\includegraphics[width=\textwidth]{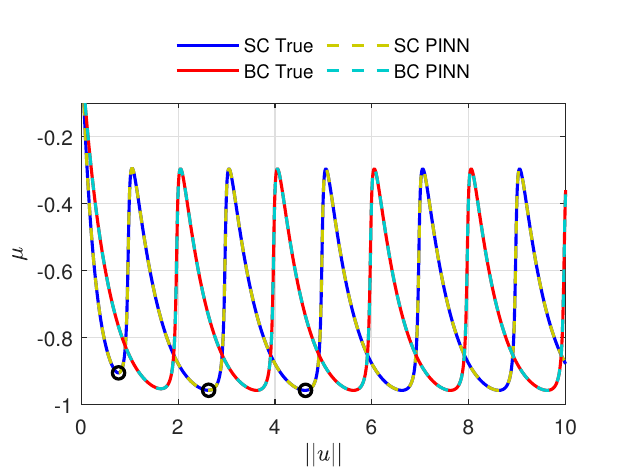}
		\caption{}
		\label{fig_1d_ac_bifurcation1}
	\end{subfigure}
	\begin{subfigure}[b]{0.49\textwidth}
		\includegraphics[width=\textwidth]{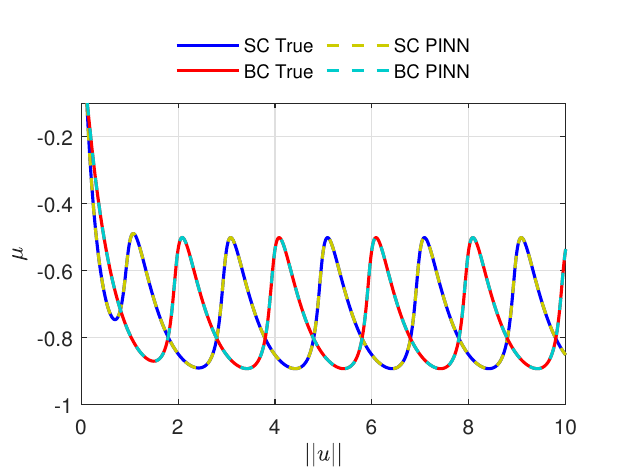}
		\caption{}
		\label{fig_1d_ac_bifurcation2}
	\end{subfigure}
	\caption{Snaking bifurcation diagrams of the one-dimensional discrete Allen--Cahn equation for (a) $c = 0.05$ and (b) $c = 0.15$. SC: site-centered, BC: bond-centered. Black circles mark the solutions shown in Figure~\ref{fig_1d_ac_solution}.}
	\label{fig_1d_ac_bifurcation}
\end{figure}

The initial solutions at $\mu = -0.1$ and $\mu = -0.12$ serve as starting points for continuation. Figure~\ref{fig_1d_ac_bifurcation} shows the resulting bifurcation diagrams for $c = 0.05$ and $c = 0.15$. The PINN-generated curves closely match the ground-truth solutions, confirming the method's ability to capture the underlying bifurcation structure accurately across both site- and bond-centered configurations.

\begin{figure}
	\centering
	
	\begin{subfigure}[b]{0.49\textwidth}
		\centering
		\includegraphics[width=\textwidth]{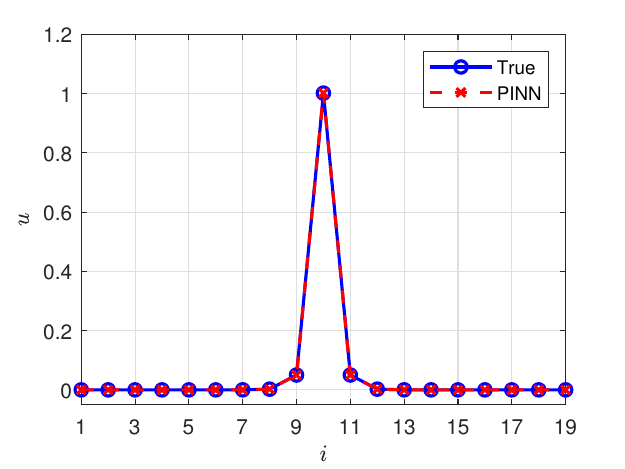}
		\caption{}
		\label{fig_1d_ac_solution1}
	\end{subfigure}
	\begin{subfigure}[b]{0.49\textwidth}
		\centering
		\includegraphics[width=\textwidth]{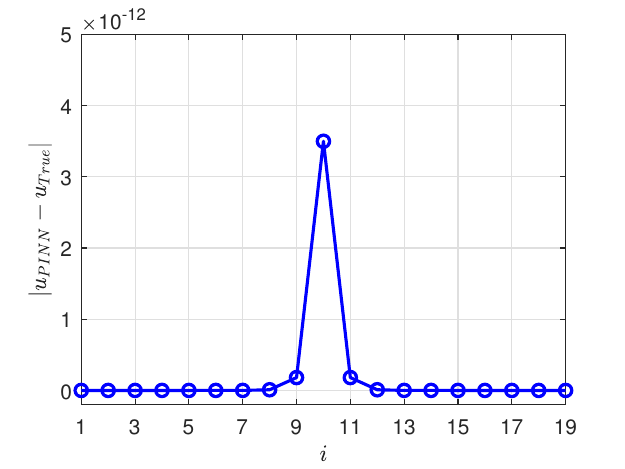}
		\caption{}
		\label{fig_1d_ac_solution2}
	\end{subfigure}
	
	\begin{subfigure}[b]{0.49\textwidth}
		\centering
		\includegraphics[width=\textwidth]{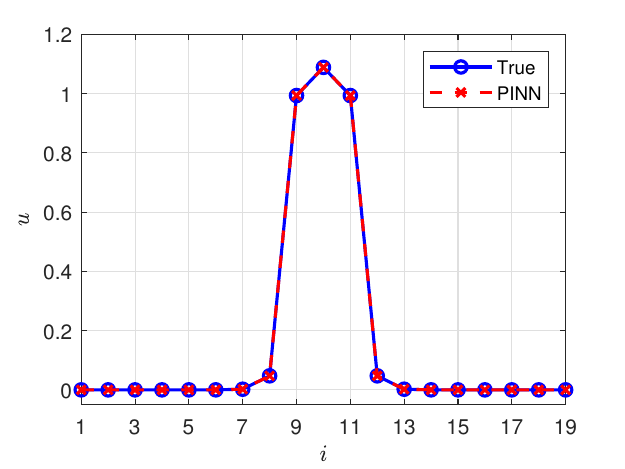}
		\caption{}
		\label{fig_1d_ac_solution3}
	\end{subfigure}
	\begin{subfigure}[b]{0.49\textwidth}
		\centering
		\includegraphics[width=\textwidth]{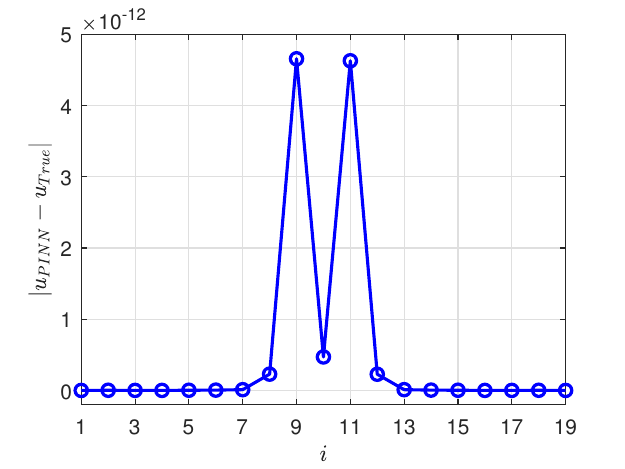}
		\caption{}
		\label{fig_1d_ac_solution4}
	\end{subfigure}
	
	\begin{subfigure}[b]{0.49\textwidth}
		\centering
		\includegraphics[width=\textwidth]{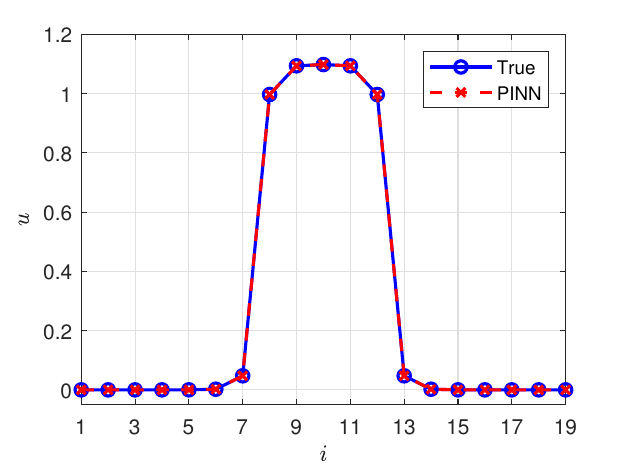}
		\caption{}
		\label{fig_1d_ac_solution5}
	\end{subfigure}
	\begin{subfigure}[b]{0.49\textwidth}
		\centering
		\includegraphics[width=\textwidth]{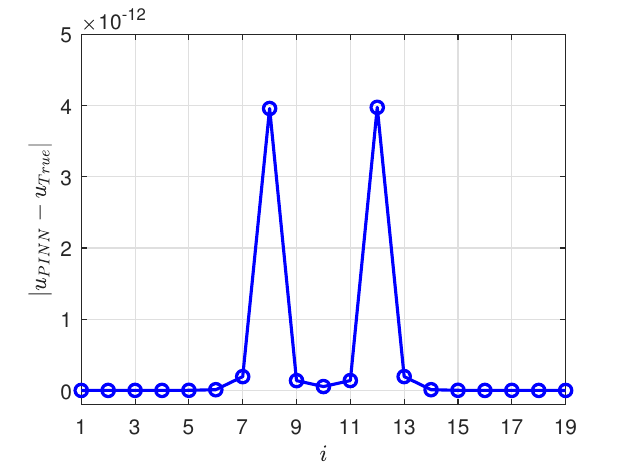}
		\caption{}
		\label{fig_1d_ac_solution6}
	\end{subfigure}
    \caption{(a), (c), (e): True and PINN solutions corresponding to the marked points in Figure~\ref{fig_1d_ac_bifurcation} for $c = 0.05$. (b), (d), (f): Absolute errors between true and PINN solutions.}
	\label{fig_1d_ac_solution}
\end{figure}

Figure~\ref{fig_1d_ac_solution} presents three representative solutions for $c = 0.05$, corresponding to the marked locations in Figure~\ref{fig_1d_ac_bifurcation}. Subfigures (a), (c), and (e) show the true and PINN solutions, while (b), (d), and (f) display the absolute differences. In all cases, the error remains below $5 \times 10^{-12}$, highlighting the exceptional accuracy of the proposed method.

\begin{figure}[t]
	\centering
	\begin{subfigure}[b]{0.49\textwidth}
		\includegraphics[width=\textwidth]{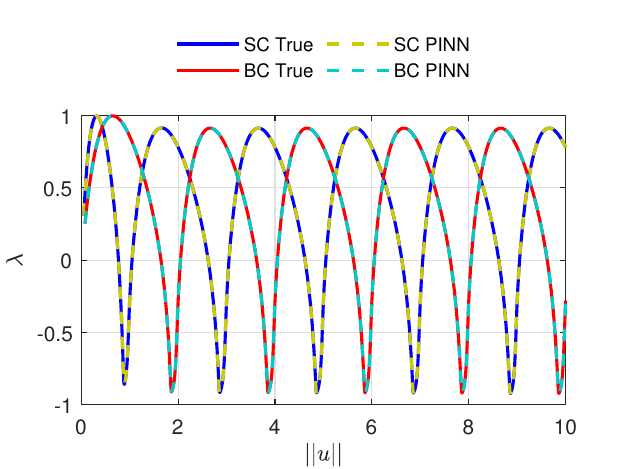}
		\caption{}
		\label{fig_1d_ac_eigenvalue1}
	\end{subfigure}
	\begin{subfigure}[b]{0.49\textwidth}
		\includegraphics[width=\textwidth]{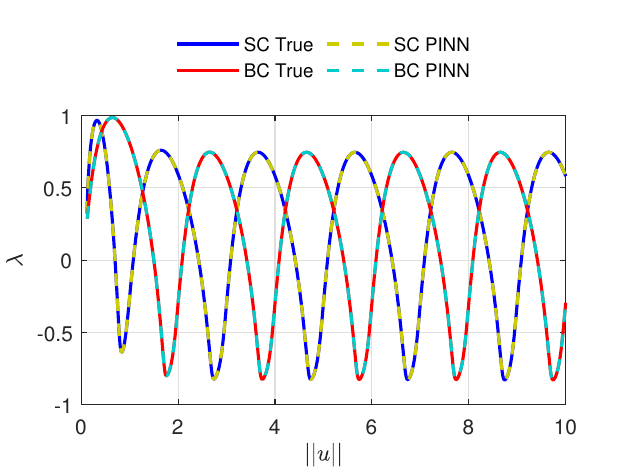}
		\caption{}
		\label{fig_1d_ac_eigenvalue2}
	\end{subfigure}
	\caption{Largest eigenvalues versus $\norm{u}$ for the one-dimensional discrete Allen--Cahn equation: (a) $c = 0.05$, (b) $c = 0.15$.}
	\label{fig_1d_ac_eigenvalue}
\end{figure}

\begin{figure}[t]
	\centering
	\begin{subfigure}[b]{0.49\textwidth}
		\includegraphics[width=\textwidth]{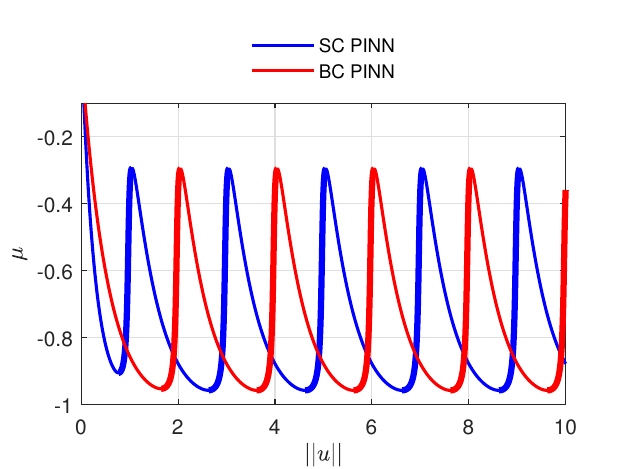}
		\caption{}
		\label{fig_1d_ac_stability1}
	\end{subfigure}
	\begin{subfigure}[b]{0.49\textwidth}
		\includegraphics[width=\textwidth]{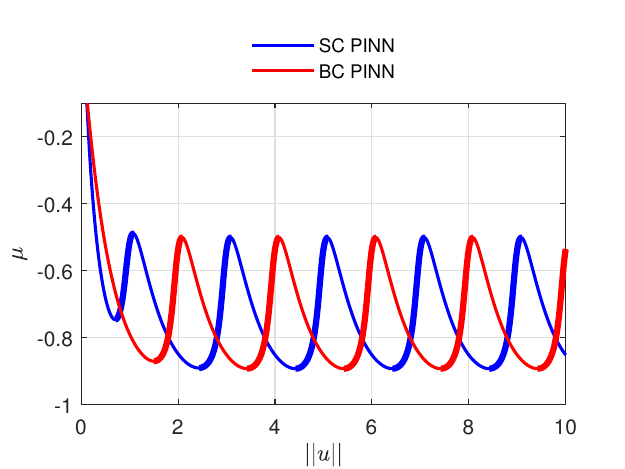}
		\caption{}
		\label{fig_1d_ac_stability2}
	\end{subfigure}
	\caption{Bifurcation diagrams with linear stability indicated: (a) $c = 0.05$, (b) $c = 0.15$. Thick lines: stable branches; thin lines: unstable branches.}
	\label{fig_1d_ac_stability}
\end{figure}

Figure~\ref{fig_1d_ac_eigenvalue} displays the largest eigenvalue as a function of $\norm{u}$ for both cases. The eigenvalues were computed using the PINN-based Sturm--Liouville approach described in Section~\ref{sec_4_3}. The transitions between stability and instability are clearly observed.
Finally, Figure~\ref{fig_1d_ac_stability} presents the stability classification of the bifurcation branches. Derived from the eigenvalue analysis in Figure~\ref{fig_1d_ac_eigenvalue}, stable and unstable regions are depicted using thick and thin lines, respectively. The alternating pattern of stability along the branches is consistent with the known behavior of snaking bifurcations.

\subsection{Two-dimensional Discrete Allen--Cahn Equation}

We now extend our analysis to the two-dimensional discrete Allen--Cahn equation, defined in Eq.~\eqref{allencahn2}. Compared to the one-dimensional case, this system exhibits more intricate and irregular snaking bifurcations, as observed in \cite{taylor2010snaking}.
For the simulation, we set $m = 8$, resulting in lattice sizes of $15 \times 15$ and $16 \times 16$ for the site- and bond-centered cases, respectively. Excluding boundary points, the corresponding Jacobian matrices for direct solution have dimensions $169 \times 169$ and $196 \times 196$.

We employ a neural network with two inputs, two hidden layers containing seven neurons each, and one output. This architecture comprises 85 trainable parameters, including biases. With PINNs, the Jacobian matrix dimensions are reduced to $169 \times 85$ and $196 \times 85$ for the site- and bond-centered cases, respectively. Moreover, in the Levenberg--Marquardt algorithm, the matrix $(J^T J + \lambda I)$ further reduces to size $85 \times 85$, regardless of configuration.

To construct the bifurcation diagrams, we use the arclength-based continuation with the following parameters: $\alpha = 10$, $\beta_1 = 1000/40$, $\beta_2 = 100$, $\gamma = 0$, and $\delta = 1$. The solution norm is defined as
\begin{equation}
	\norm{u} = \frac{1}{1+\sqrt{1+\mu}} \sum_{i=1}^{n} \sum_{j=1}^{n} u_{i,j}^2,
\end{equation}
and the diagrams are computed up to $\norm{u} = 40$ for both solution types. Sensitivity analysis for $\beta_1$ and $\beta_2$ is provided in Appendix~\ref{appendix_beta}.

\begin{figure}
	\centering
	\begin{subfigure}[b]{0.49\textwidth}
		\includegraphics[width=\textwidth]{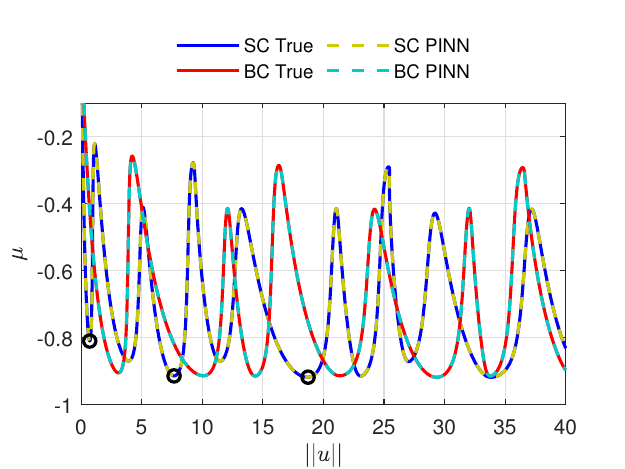}
		\caption{}
		\label{fig_2d_ac_bifurcation1}
	\end{subfigure}
	\begin{subfigure}[b]{0.49\textwidth}
		\includegraphics[width=\textwidth]{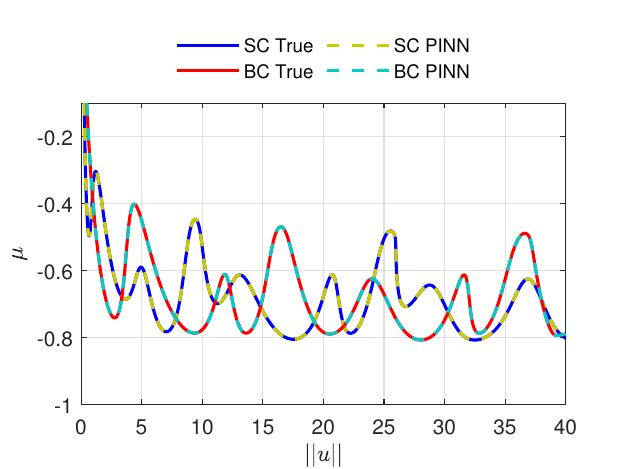}
		\caption{}
		\label{fig_2d_ac_bifurcation2}
	\end{subfigure}
	\caption{Snaking bifurcation diagrams of the two-dimensional discrete Allen--Cahn equation for (a) $c = 0.05$ and (b) $c = 0.15$. SC: site-centered, BC: bond-centered. Black circles mark the solutions shown in Figure~\ref{fig_2d_ac_solution}.}
	\label{fig_2d_ac_bifurcation}
\end{figure}

The diagrams are initiated using two solutions at $\mu = -0.1$ and $\mu = -0.1 - 1/\beta_2$. These serve as starting points for the pseudo-arclength continuation. Figure~\ref{fig_2d_ac_bifurcation} presents the resulting bifurcation diagrams for $c = 0.05$ and $c = 0.15$. In both site- and bond-centered cases, the PINN-generated results are nearly indistinguishable from those of the direct method, demonstrating the model's effectiveness. Unlike the one-dimensional case, the snaking patterns in two dimensions are more complex and irregular, underscoring the challenge of accurate bifurcation tracking in higher dimensions.

\begin{figure}
	\centering
	
	\begin{subfigure}[b]{0.49\textwidth}
		\centering
		\includegraphics[width=\textwidth]{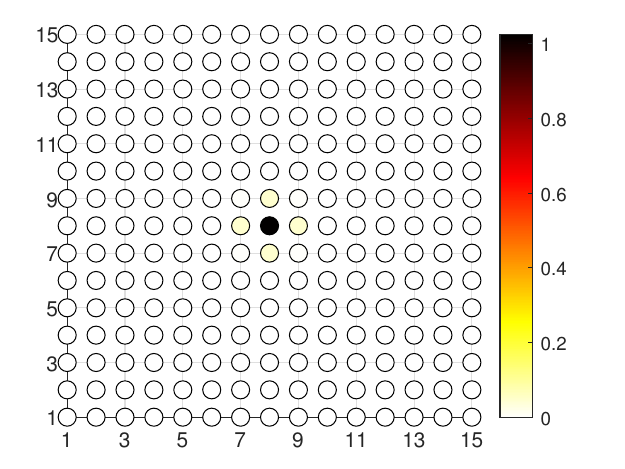}
		\caption{}
		\label{fig_2d_ac_solution1}
	\end{subfigure}
	\begin{subfigure}[b]{0.49\textwidth}
		\centering
		\includegraphics[width=\textwidth]{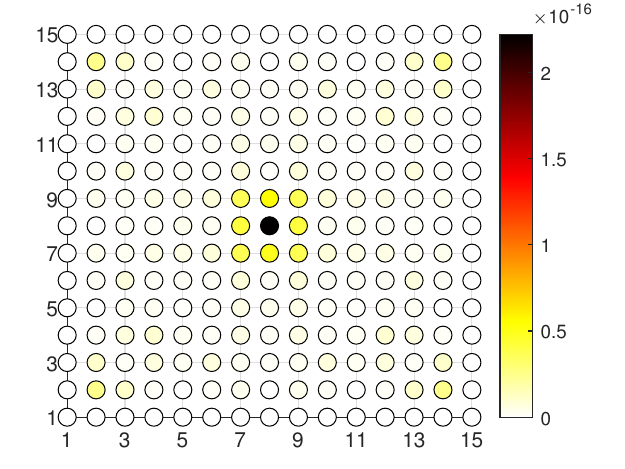}
		\caption{}
		\label{fig_2d_ac_solution2}
	\end{subfigure}
	
	\begin{subfigure}[b]{0.49\textwidth}
		\centering
		\includegraphics[width=\textwidth]{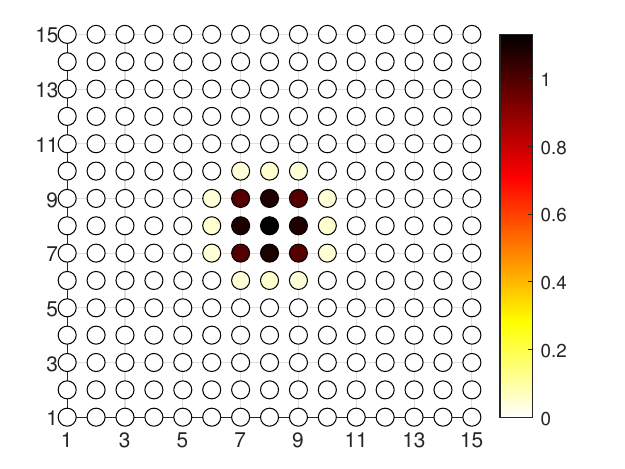}
		\caption{}
		\label{fig_2d_ac_solution3}
	\end{subfigure}
	\begin{subfigure}[b]{0.49\textwidth}
		\centering
		\includegraphics[width=\textwidth]{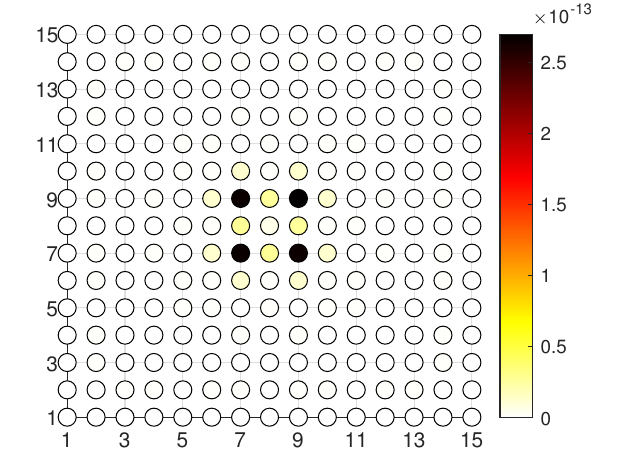}
		\caption{}
		\label{fig_2d_ac_solution4}
	\end{subfigure}
	
	\begin{subfigure}[b]{0.49\textwidth}
		\centering
		\includegraphics[width=\textwidth]{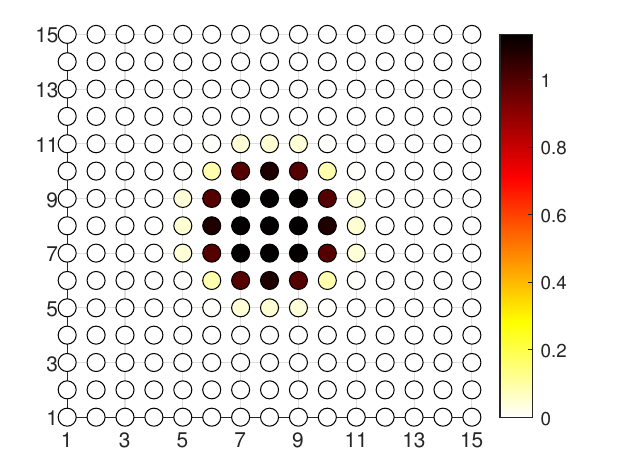}
		\caption{}
		\label{fig_2d_ac_solution5}
	\end{subfigure}
	\begin{subfigure}[b]{0.49\textwidth}
		\centering
		\includegraphics[width=\textwidth]{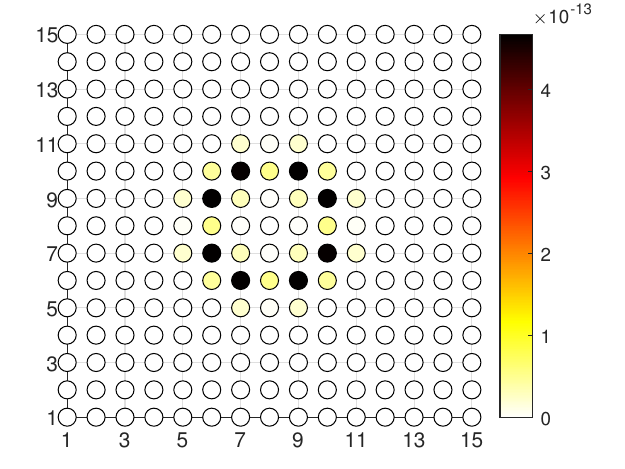}
		\caption{}
		\label{fig_2d_ac_solution6}
	\end{subfigure}
		\caption{(a), (c), (e): PINN solutions corresponding to black circles in Figure~\ref{fig_2d_ac_bifurcation} for $c = 0.05$. (b), (d), (f): Absolute difference between true and PINN solutions.}
	\label{fig_2d_ac_solution}
\end{figure}

\begin{figure}[t]
	\centering
	\begin{subfigure}[b]{0.49\textwidth}
		\includegraphics[width=\textwidth]{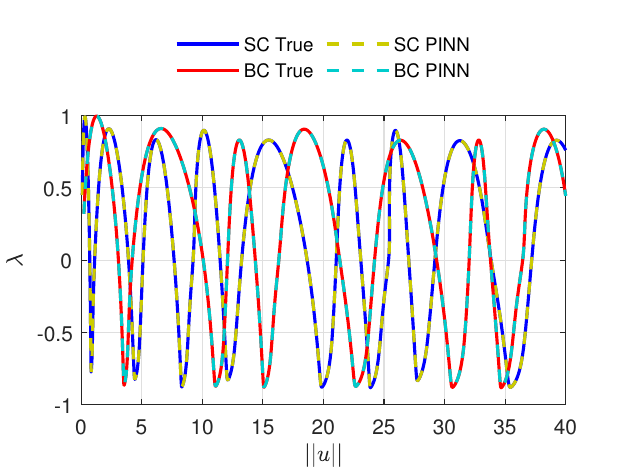}
		\caption{}
		\label{fig_2d_ac_eigenvalue1}
	\end{subfigure}
	\begin{subfigure}[b]{0.49\textwidth}
		\includegraphics[width=\textwidth]{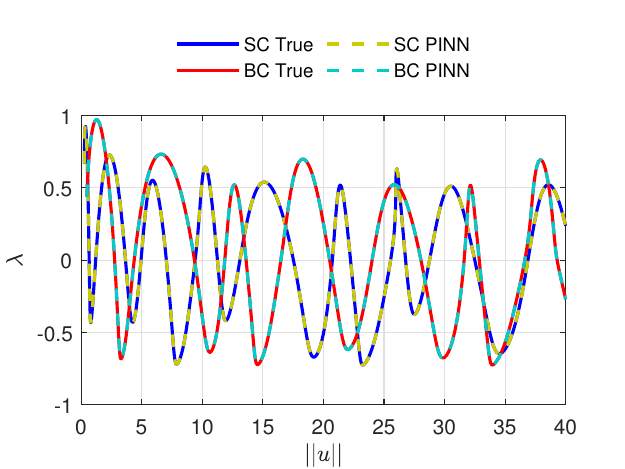}
		\caption{}
		\label{fig_2d_ac_eigenvalue2}
	\end{subfigure}
	\caption{Largest eigenvalue versus $\norm{u}$ for the two-dimensional discrete Allen--Cahn equation: (a) $c = 0.05$, (b) $c = 0.15$. SC: site-centered, BC: bond-centered.}
	\label{fig_2d_ac_eigenvalue}
\end{figure}

\begin{figure}[H]
	\centering
	\begin{subfigure}[b]{0.49\textwidth}
		\includegraphics[width=\textwidth]{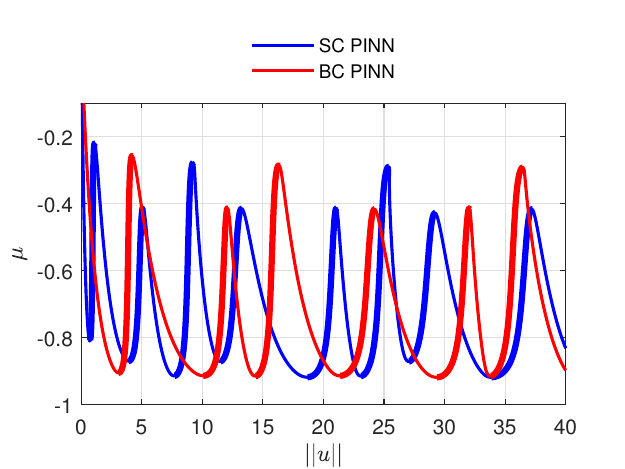}
		\caption{}
		\label{fig_2d_ac_stability1}
	\end{subfigure}
	\begin{subfigure}[b]{0.49\textwidth}
		\includegraphics[width=\textwidth]{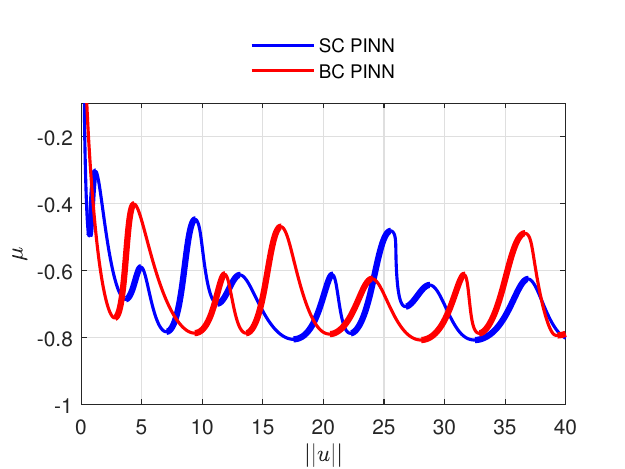}
		\caption{}
		\label{fig_2d_ac_stability2}
	\end{subfigure}
	\caption{Snaking bifurcation diagrams with linear stability indicated for (a) $c = 0.05$ and (b) $c = 0.15$. SC: site-centered, BC: bond-centered. Thick lines: stable branches; thin lines: unstable branches.}
	\label{fig_2d_ac_stability}
\end{figure}

Figure~\ref{fig_2d_ac_solution} shows three representative solutions at selected values of $\norm{u}$ for $c = 0.05$, marked in Figure~\ref{fig_2d_ac_bifurcation}. The left column displays the PINN approximations, while the right column presents the absolute errors compared to the true solutions. In all instances, the absolute error remains below $5 \times 10^{-13}$, confirming the exceptional precision of the PINN method.

Figure~\ref{fig_2d_ac_eigenvalue} displays the largest eigenvalue as a function of $\norm{u}$, computed using the PINN-based Sturm--Liouville approach. As expected, the eigenvalues fluctuate with the bifurcation branches, reflecting the stability characteristics of the solutions.
Finally, Figure~\ref{fig_2d_ac_stability} shows the linear stability classification derived from the eigenvalue analysis in Figure~\ref{fig_2d_ac_eigenvalue}. As in the one-dimensional case, stability alternates along the bifurcation curve, with thick lines denoting stable branches and thin lines indicating unstable ones. The results reaffirm the PINN framework’s ability to accurately resolve both the bifurcation structure and the associated stability transitions, even in more complex two-dimensional systems.

\subsection{Three-dimensional Discrete Allen--Cahn Equation}

We now extend our study to higher-dimensional versions of the discrete Allen--Cahn equation, beginning with the three-dimensional case. Due to the computational cost of constructing complete bifurcation diagrams in higher dimensions, we restrict our analysis to a single parameter value in the site-centered case.

We consider the three-dimensional problem at a fixed parameter $\mu = -0.5$ with $c = 0.05$. Choosing $m = 8$ results in a lattice with $15^3$ points. Excluding the boundary points, the system consists of $13^3 = 2197$ nonlinear equations and variables. Directly solving the system requires computing a Jacobian matrix of size $2197 \times 2197$.

\begin{figure}[h]
	\centering
	\begin{subfigure}[b]{0.49\textwidth}
		\includegraphics[width=\textwidth]{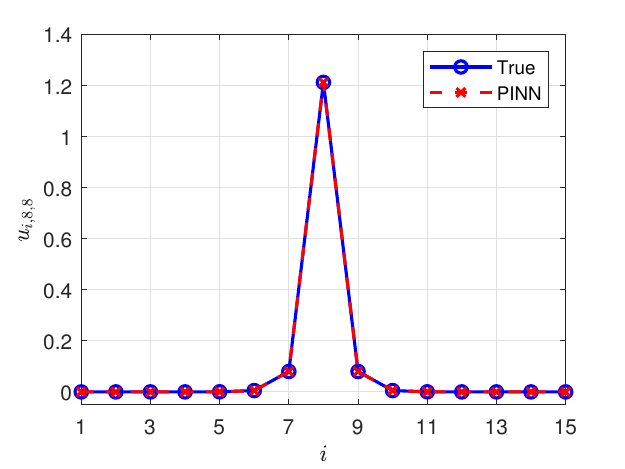}
		\caption{}
		\label{fig_3d_ac_solution}
	\end{subfigure}
	\begin{subfigure}[b]{0.49\textwidth}
		\includegraphics[width=\textwidth]{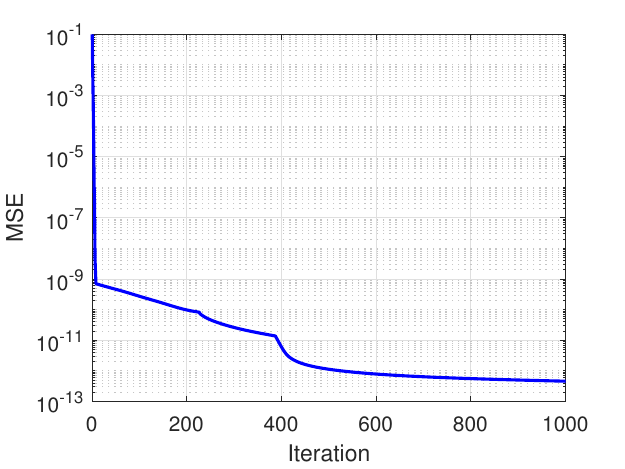}
		\caption{}
		\label{fig_3d_ac_mse}
	\end{subfigure}
	\caption{(a) True and PINN solutions for the three-dimensional discrete Allen--Cahn equation, plotted as $u_{i,8,8}$ for $i = 1, \dots, 15$. (b) Mean squared error (MSE) during training.}
	\label{fig_3d_ac}
\end{figure}

To approximate the solution using PINNs, we employ a neural network with three inputs, two hidden layers each containing 10 neurons, and one output, resulting in 161 trainable parameters including biases. The corresponding Jacobian matrix has dimensions $2197 \times 161$, significantly smaller than in the direct method.
Figure~\ref{fig_3d_ac_solution} shows a one-dimensional slice of the solution, $u_{i,8,8}$, for visualization. The training convergence, shown in Figure~\ref{fig_3d_ac_mse}, confirms that the PINN accurately solves the three-dimensional problem, achieving an MSE on the order of $10^{-13}$.

\subsection{Four-dimensional Discrete Allen--Cahn Equation}

We next consider the four-dimensional case at fixed parameter $\mu = -0.5$ and $c = 0.05$. Choosing $m = 8$ yields a lattice with $15^4$ points and $13^4 = 28561$ internal equations and variables. Solving the system directly would require manipulating a $28561 \times 28561$ Jacobian, which exceeds typical memory capacities.

\begin{figure}[t]
	\centering
	\begin{subfigure}[b]{0.49\textwidth}
		\includegraphics[width=\textwidth]{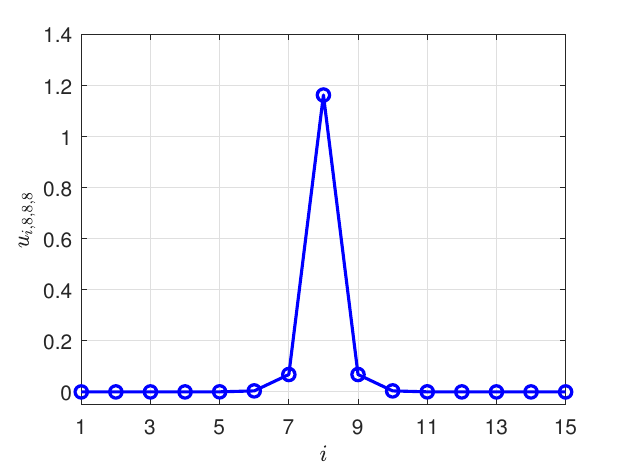}
		\caption{}
		\label{fig_4d_ac_solution}
	\end{subfigure}
	\begin{subfigure}[b]{0.49\textwidth}
		\includegraphics[width=\textwidth]{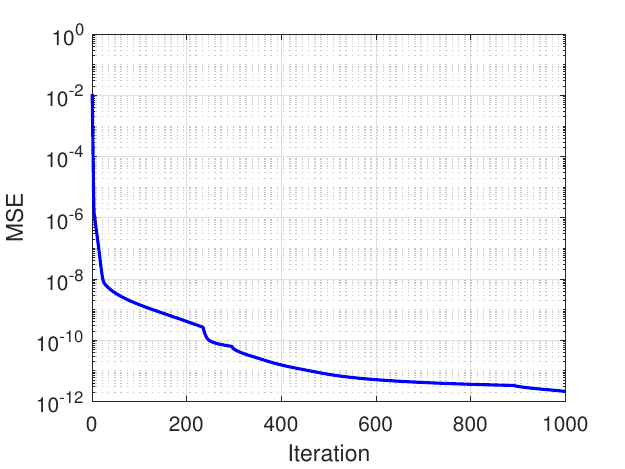}
		\caption{}
		\label{fig_4d_ac_mse}
	\end{subfigure}
	\caption{(a) PINN solution for the four-dimensional Allen--Cahn equation shown as $u_{i,8,8,8}$. (b) MSE during training.}
	\label{fig_4d_ac}
\end{figure}

Using a PINN with four inputs, two hidden layers of 10 neurons each, and one output, the number of weights is reduced to 171. The resulting Jacobian matrix has dimensions $28561 \times 171$, a significant reduction that makes the computation feasible.
Figure~\ref{fig_4d_ac_solution} presents a one-dimensional slice of the solution, while Figure~\ref{fig_4d_ac_mse} shows that the MSE reaches the order of $10^{-12}$, confirming the high accuracy of the PINN approach in four dimensions.

\subsection{Five-dimensional Discrete Allen--Cahn Equation}

Finally, we address the five-dimensional case using the same fixed parameters: $\mu = -0.5$ and $c = 0.05$. With $m = 8$, the site-centered lattice includes $15^5$ points and $13^5 = 371293$ internal variables and equations.

We use a PINN with five inputs, two hidden layers of 10 neurons each, and one output, resulting in 181 trainable parameters. The Jacobian matrix would nominally be $371293 \times 181$, but due to its size, we adopt the stochastic approach described in Section~\ref{stochastic_newton}. In each iteration, we randomly sample 1000 equations and always include the central function at $\textbf{\textit{i}} = (8, 8, 8, 8, 8)$ to increase accuracy.
This reduces the Jacobian size to $1001 \times 181$, and the Levenberg--Marquardt update involves matrices of size $181 \times 181$ and $181 \times 1$, substantially improving computational efficiency.

\begin{figure}[t]
	\centering
	\begin{subfigure}[b]{0.49\textwidth}
		\includegraphics[width=\textwidth]{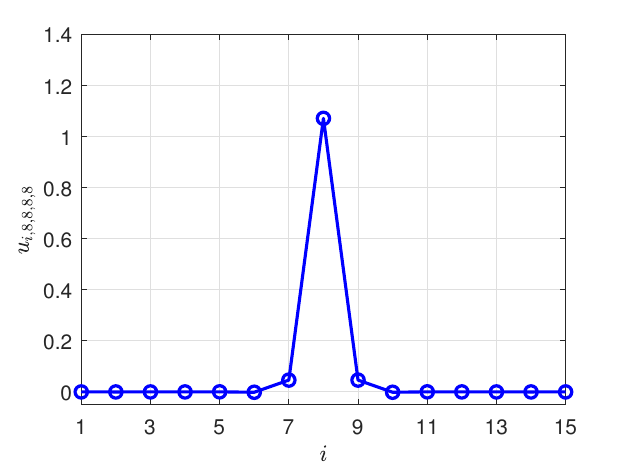}
		\caption{}
		\label{fig_5d_ac_solution}
	\end{subfigure}
	\begin{subfigure}[b]{0.49\textwidth}
		\includegraphics[width=\textwidth]{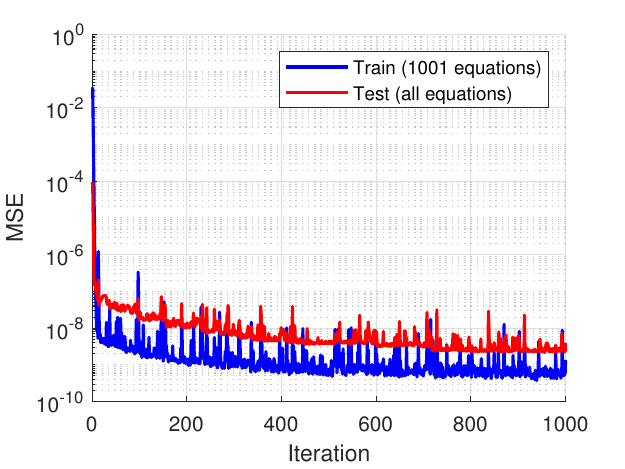}
		\caption{}
		\label{fig_5d_ac_mse}
	\end{subfigure}
	\caption{(a) PINN solution for the five-dimensional Allen--Cahn equation shown as $u_{i,8,8,8,8}$. (b) MSE during training: blue for training set (1001 equations), red for full system (371293 equations).}
	\label{fig_5d_ac}
\end{figure}

Figure~\ref{fig_5d_ac_solution} shows the PINN solution slice, and Figure~\ref{fig_5d_ac_mse} presents the MSE for both training and test sets. Despite the stochastic updates, the final test MSE is approximately $10^{-9}$, indicating that the PINN yields a sufficiently accurate solution. The non-monotonic MSE behavior arises from the stochastic nature of training, unlike the smooth convergence seen in lower dimensions where full Jacobians are used. Furthermore, the effect of varying the number of equations used during training is reported in Table~\ref{tab_5d}.

\begin{table}[t]
	\centering
	\caption{Effect of training set size s$(S_k)$ on test MSE (full system) and computation time for the five-dimensional case over 1000 iterations.}
	\label{tab_5d}
	\begin{tabular}{|c|c|c|c|}
		\hline
		s$(S_k)$ & $n_{\text{iter}}$ & Test MSE & Time (s) \\ \hline
		101   & 1000 & $4.1 \times 10^{-8}$ & 88 \\
		1001  & 1000 & $4.0 \times 10^{-9}$ & 108 \\
		10001 & 1000 & $4.4 \times 10^{-10}$ & 301 \\ \hline
	\end{tabular}
\end{table}

\begin{table}[t]
	\centering
	\caption{Summary of PINN architectures from one to five dimensions.}
	\label{tab_1d_5d}
	\begin{tabular}{|c|c|c|c|c|c|c|c|}
		\hline
		$d$ & PINN & $n_W$ & s$(J)$ & s$(J^T J+\lambda I)$ & $n_{\text{iter}}$ & MSE & Time (s) \\ \hline
		1 & (1,4,4,1) & 33  & $17 \times 33$     & $33 \times 33$     & 1000 & $2.1 \times 10^{-32}$ & 0.1 \\
		2 & (2,7,7,1) & 85  & $169 \times 85$    & $85 \times 85$     & 1000 & $1.0 \times 10^{-13}$ & 10  \\
		3 & (3,10,10,1) & 161 & $2197 \times 161$ & $161 \times 161$   & 1000 & $4.6 \times 10^{-13}$ & 33  \\
		4 & (4,10,10,1) & 171 & $28561 \times 171$ & $171 \times 171$  & 1000 & $2.1 \times 10^{-12}$ & 552 \\
		5 & (5,10,10,1) & 181 & $1001 \times 181$ & $181 \times 181$   & 1000 & $4.0 \times 10^{-9}$  & 108 \\ \hline
	\end{tabular}
\end{table}

\begin{table}[t]
	\centering
	\caption{Impact of hidden layer size on MSE and training time for the three-dimensional case.}
	\label{tab_3D}
	\begin{tabular}{|c|c|c|c|c|c|c|c|}
		\hline
		$d$ & PINN & $n_W$ & s$(J)$ & s$(J^T J + \lambda I)$ & $n_{\text{iter}}$ & MSE & Time (s) \\ \hline
		3 & (3,5,5,1) & 56  & $2197 \times 56$  & $56 \times 56$   & 1000 & $1.1 \times 10^{-12}$ & 11  \\
		3 & (3,10,5,1) & 101 & $2197 \times 101$ & $101 \times 101$ & 1000 & $1.1 \times 10^{-12}$ & 22  \\
		3 & (3,10,10,1) & 161 & $2197 \times 161$ & $161 \times 161$ & 1000 & $4.6 \times 10^{-13}$ & 33  \\
		3 & (3,15,10,1) & 231 & $2197 \times 231$ & $231 \times 231$ & 1000 & $8.5 \times 10^{-13}$ & 53  \\
		3 & (3,15,15,1) & 316 & $2197 \times 316$ & $316 \times 316$ & 1000 & $2.2 \times 10^{-13}$ & 81  \\
		3 & (3,20,15,1) & 411 & $2197 \times 411$ & $411 \times 411$ & 1000 & $9.0 \times 10^{-13}$ & 128 \\
		3 & (3,20,20,1) & 521 & $2197 \times 521$ & $521 \times 521$ & 1000 & $1.9 \times 10^{-12}$ & 210 \\ \hline
	\end{tabular}
\end{table}

Tables~\ref{tab_1d_5d} and \ref{tab_3D} summarize the architectural, numerical, and performance aspects of the PINNs across dimensions. The results indicate that the proposed framework scales well with dimensionality and provides accurate solutions, even in five-dimensional systems where traditional methods are computationally prohibitive.

\section{Conclusion}\label{sec6}

In this work, we proposed the use of physics-informed neural networks (PINNs) to address nonlinear systems arising from discrete nonlinear lattices. The scope of the study encompassed three core tasks: obtaining steady-state solutions, constructing bifurcation diagrams, and performing linear stability analysis. To that end, we formulated a PINN-based framework that processes lattice-based systems as input and solves the resulting nonlinear equations using the Levenberg--Marquardt algorithm. For high-dimensional settings, where conventional methods struggle due to memory constraints and computational cost, we introduced a stochastic optimization strategy to improve scalability and efficiency.

In the one- and two-dimensional cases, we demonstrated the efficacy of PINNs in capturing the snaking bifurcation structure and computing the associated linear stability. The continuation method was implemented via pseudo-arclength continuation, which was adapted to the PINN framework through the inclusion of an auxiliary constraint. This allowed for a sequence of neural network weights to be computed iteratively, thereby enabling the construction of full bifurcation diagrams. Depending on the nature of the solution curve, the method was flexibly applied in two variants—either assuming or not assuming a bijective relation between the norm and the solution.

For stability analysis, we proposed a PINN formulation to compute the principal eigenvalue of the linearized system. By enforcing positivity in the network output, we leveraged Sturm--Liouville theory to isolate the largest eigenvalue associated with a non-sign-changing eigenvector. Numerical results showed excellent agreement between the PINN predictions and reference solutions, affirming the robustness and accuracy of the approach.

In three-, four-, and five-dimensional settings, we extended the framework by employing compact neural networks with only modest growth in the number of trainable parameters, despite the exponential increase in system size. Each additional spatial dimension was accommodated by augmenting the input layer, while the number of hidden neurons remained fixed. This design choice helped preserve computational efficiency and kept the size of the Jacobian matrix tractable. Our experiments confirmed that even in high-dimensional spaces, PINNs can produce accurate solutions with considerably lower memory usage than traditional methods. In the five-dimensional case, we further applied a stochastic Levenberg--Marquardt algorithm, evaluating only a subset of equations per iteration, which significantly accelerated training and reduced computational overhead.

Looking ahead, several promising avenues can be explored to extend this work. One direction involves applying the proposed PINN methodology to other nonlinear lattice models that exhibit complex bifurcation phenomena. Notable examples include the discrete Swift--Hohenberg equation \cite{kusdiantara2017homoclinic, burke2006localized} and the discrete nonlinear Schr\"odinger equation, which models quantum droplets and bubbles \cite{kusdiantara2024analysis}. Similarly, adapting the framework to non-standard lattice geometries, such as Lieb \cite{kusdiantara2022snakes}, honeycomb, and triangular lattices \cite{kusdiantara2019snakes}, offers further opportunities for exploration. 

Beyond lattice-based models, the techniques developed here may be generalized to a broader class of nonlinear systems where high dimensionality poses a challenge. The inherent flexibility of neural network architectures makes them particularly well-suited for capturing solution manifolds in such settings. Furthermore, in problems where analytical solutions are known, as in the integrable Ablowitz--Ladik lattice \cite{zhu2022neural}, PINNs may serve not only as numerical solvers but also as tools for symbolic discovery. Recent studies have demonstrated that neural networks can be trained to approximate and even recover closed-form solutions \cite{dana2025approximate, mazraeh2025innovative, mazraeh2024gepinn, mazraeh2025three}.

In summary, this work highlights the potential of PINNs as an effective and scalable tool for solving nonlinear systems on lattices, particularly in high-dimensional and bifurcation-sensitive regimes. The combination of accuracy, adaptability, and computational efficiency makes them a compelling alternative to traditional numerical methods, especially for problems where conventional solvers become impractical.

\section*{CRediT authorship contribution statement}
\textbf{Muhammad Luthfi Shahab:} Conceptualization, Formal Analysis, Software, Visualization, Writing - original draft, 
\textbf{Fidya Almira Suheri:} Conceptualization, Software,
\textbf{Rudy Kusdiantara:} Conceptualization, Writing - review \& editing,
\textbf{Hadi Susanto:} Conceptualization, Supervision, Writing - review \& editing.


\section*{Declaration of competing interest} 
The authors declare that they have no known competing financial interests or personal relationships that could have appeared to influence the work reported in this paper.

\section*{Declaration of generative AI and AI-assisted technologies in the writing process}

During the preparation of this work the authors used Grammarly and ChatGPT in order to improve language and readability. After using these tools/services, the authors reviewed and edited the content as needed and took full responsibility for the content of the publication.


\section*{Acknowledgement}
MLS is supported by a four-year Doctoral Research and Teaching Scholarship (DRTS) from Khalifa University. RK acknowledges Riset Utama PPMI FMIPA 2024 (617I/IT1.C02/KU/2024). HS acknowledged support by Khalifa University through a Competitive Internal Research Awards Grant (No.\ 8474000413/CIRA-2021-065) and Research \& Innovation Grants (No.\ 8474000617/RIG-S-2023-031 and No.\ 8474000789/RIG-S-2024-070).

\bibliographystyle{elsarticle-num} 
\bibliography{main_references}

\appendix

\section{Effect of Varying $\alpha$ and $\gamma$} \label{appendix_alpha_gamma}

\begin{figure}[h]
	\centering
	
	\begin{subfigure}[b]{0.49\textwidth}
		\centering
		\includegraphics[width=\textwidth]{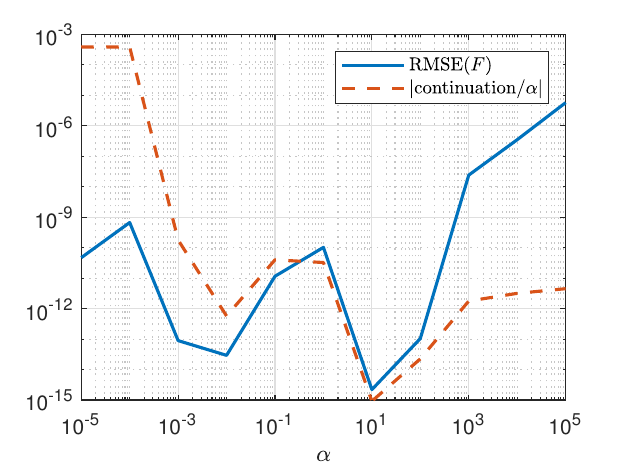}
		\caption{}
		\label{fig_1d_ac_alpha}
	\end{subfigure}
	\begin{subfigure}[b]{0.49\textwidth}
		\centering
		\includegraphics[width=\textwidth]{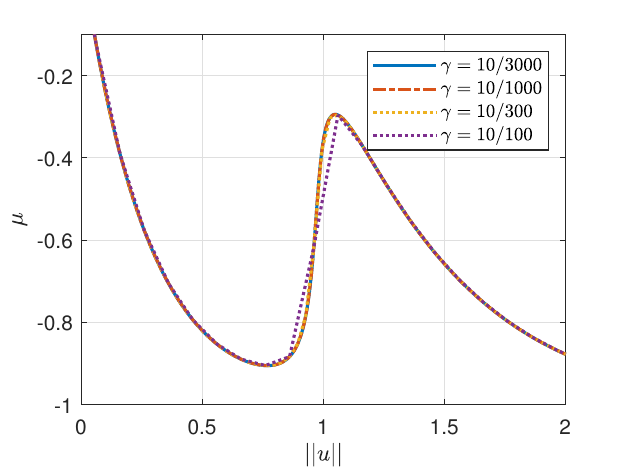}
		\caption{}
		\label{fig_1d_ac_gamma}
	\end{subfigure}
	
	\caption{(a) Influence of $\alpha$ and (b) $\gamma$ on the one-dimensional Allen-Cahn equation.}
	\label{fig_1d_ac_alpha_gamma}
\end{figure}

The parameter $\alpha$ governs the relative weighting of the continuation equation in Eq.\ \eqref{arclength1} with respect to the original nonlinear system $F(u,\mu) = 0$. When $\alpha$ is too small, the continuation constraint is weak and may be violated; when too large, it compromises the accuracy of the solution to the original system. Figure~\ref{fig_1d_ac_alpha} illustrates the effect of varying $\alpha$ across several orders of magnitude, from $10^{-5}$ to $10^{5}$. The results correspond to the site-centered case with $c = 0.05$ in the one-dimensional Allen-Cahn model. Here, the RMSE of the nonlinear system (blue line) and the scaled error in the continuation equation, $\abs{\text{continuation} / \alpha}$, are plotted. The value $\alpha = 10$ achieves a desirable trade-off between these competing errors.

The parameter $\gamma$ defines the step size in the norm of the solution between successive continuation steps. For example, setting $\gamma = 10/1000$ implies roughly 1000 steps to reach $\norm{u} = 10$. A small $\gamma$ improves resolution near sharp turning points but increases computational cost, while a large $\gamma$ may compromise accuracy. Figure~\ref{fig_1d_ac_gamma} demonstrates the consequences of varying $\gamma$ for the one-dimensional Allen-Cahn equation, where bifurcation diagrams are computed up to $\norm{u} = 2$. The diagram for $\gamma = 10/100$ fails to resolve the second turning point accurately, and even $\gamma = 10/300$ exhibits reduced smoothness. In contrast, both $\gamma = 10/1000$ and $\gamma = 10/3000$ provide smooth and consistent results. Computational times (in seconds) for each setting are 4.3250, 4.4447, 1.1761, and 9.0523, respectively, with $\gamma = 10/1000$ offering the most efficient balance between accuracy and speed.

\section{Effect of Varying $\beta_1$ and $\beta_2$} \label{appendix_beta}

\begin{figure}[b]
	\centering
	\includegraphics[width=\textwidth]{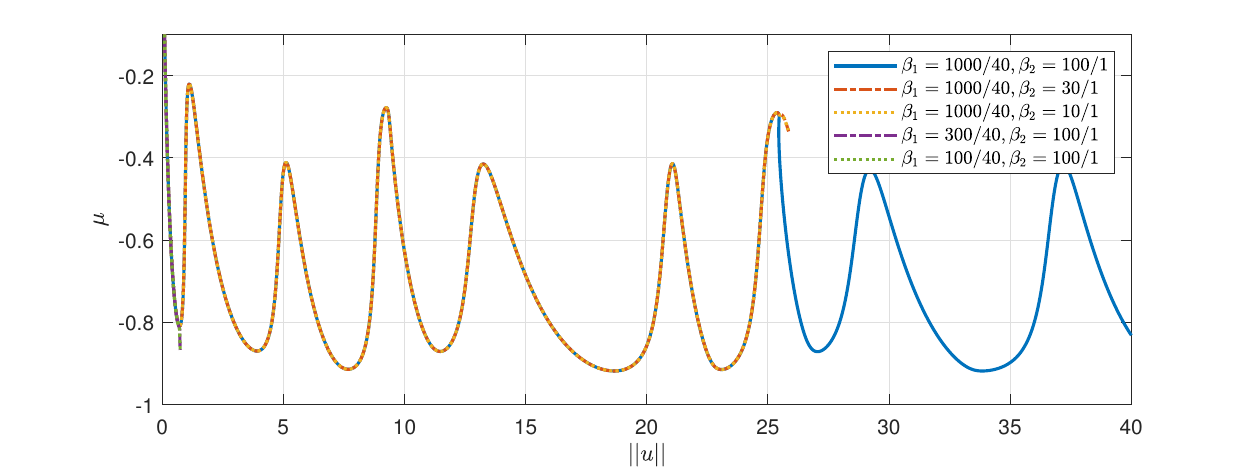}
	\caption{Effect of $\beta_1$ and $\beta_2$ on the two-dimensional Allen-Cahn equation.}
	\label{fig_2d_ac_beta}
\end{figure}

The parameters $\beta_1$ and $\beta_2$ determine the relative step sizes in the norm of the solution and the parameter $\mu$, respectively. In the two-dimensional case, we use $\beta_1 = 1000/40$ and $\beta_2 = 100/1$, which correspond to approximately 1000 steps to reach $\norm{u} = 40$ and 100 steps to traverse one unit in $\mu$. Figure~\ref{fig_2d_ac_beta} shows the impact of varying these parameters for the site-centered case with $c = 0.05$.
Smaller values of $\beta_1$ or $\beta_2$ reduce the number of continuation steps but may fail to resolve sharp turning points accurately. This is evident in the figure, where insufficient resolution leads to distorted or incomplete bifurcation structures. Conversely, increasing $\beta_1$ or $\beta_2$ improves accuracy at the cost of additional computation.
It is also worth noting that the parameter $\delta$ in Eq.\ \eqref{arclength_NN} functions in tandem with $\beta_1$ and $\beta_2$ to control the arclength step. In practice, varying $\delta$ independently is unnecessary, as its effect is implicitly governed by appropriate choices of $\beta_1$ and $\beta_2$.

\end{document}